\documentclass[10pt,twoside]{article}
\usepackage{graphicx}
\usepackage{amsmath}
\usepackage{titlesec}
\usepackage{amscd}
\usepackage{amssymb}
\usepackage{latexsym}
\usepackage{pgf,pgfarrows}

\newcommand{\R}{\mathbb R}
\newcommand{\C}{\mathbb C}
\renewcommand\thesection{\arabic{section}}
\titleformat{\section}{\normalfont\Large\bfseries}{\thesection.}{0.5em}{}
\newenvironment{proof}
{\noindent {\em Proof. }} {\hfill $\Box$}

\newtheorem{thm}{Theorem}[section]
\numberwithin{equation}{section}
\newtheorem{prop}{Proposition}[section]
\newtheorem{defn}{Definition}[section]
\newtheorem{lem}{Lemma}[section]
\title{\bf Lectures on mean curvature flows in higher codimensions}

\author{Mu-Tao Wang\thanks{Department of Mathematics \& Columbia
University, 2990 Broadway, New York, NY 10027, USA. E-mail: mtwang@math.columbia.edu}}
\date{May 26, 2008}

\begin{document}
\maketitle

\begin{abstract}\vskip 3mm\footnotesize
\noindent Mean curvature flows of hypersurfaces have been extensively studied and there are various different approaches and many beautiful results. However, relatively little is known about mean curvature flows of submanifolds of higher codimensions. This notes starts with some basic materials on submanifold geometry, and then introduces mean curvature flows in general dimensions and co-dimensions. The related techniques in the so called ``blow-up" analysis are also discussed. At the end, we present some global existence and convergence results for mean curvature flows of two-dimensional surfaces in four-dimensional ambient spaces.

\vskip 4.5mm


\end{abstract}

\vskip 12mm

\section{Basic materials}
\subsection{Connections, curvature, and the Laplacian}
A Riemannian manifold is a differentiable manifold equipped with a smooth inner product on the tangent bundle. Suppose $M$ is a Riemannian manifold of dimensional $N$ with a Riemannian metric $\langle\cdot, \cdot \rangle$.

There is a unique Levi-Civita connection $\nabla$ on the tangent bundle of $M$ that is compatible with the differentiable structure, i.e.
\[\nabla_X Y-\nabla_Y X=[X, Y]=XY-YX\] and with the metric, i.e.
\[X\langle Y, Z\rangle=\langle \nabla_X Y,Z\rangle+\langle Y,
\nabla_X Z\rangle.\]

We can extend this connection to the sections of any tensor bundles by requiring the Leibnitz rule and the compatibility with contractions
of tensors.

Suppose $X, Y, Z, W, U$, and $V$ are tangent vector fields on $M$. The
Riemannian curvature tensor is defined to be
\[R(X, Y) Z\equiv-\nabla_X\nabla_Y Z+\nabla_Y\nabla_X Z+\nabla_{[X, Y]} Z.\]

$R(X, Y, Z, W)\equiv\langle R(X, Y) Z, W\rangle$ has the following symmetries:
\[ R(X, Y, Z, W)=-R(Y, X, Z, W) \,\,\text{and}\,\,R(X, Y, Z, W)=R(Z, W, X, Y).\]
The first Bianchi identity:\[R(X, Y, Z, W)+R(X, Z, W, Y)+R(X, W, Y, Z)=0.\]
The second Bianchi identity: \[(\nabla_X R)(U, V, Y, Z)+(\nabla_Y R)(U, V, Z, X)+(\nabla_Z R)(U, V, X, Y)=0.\]

Given a smooth function $f$ on $M$, the gradient of $f$, $\nabla f$, is a tangent vector field that satisfies $\langle \nabla f, X\rangle=df(X)$ for any $X$. The Hessian of $f$, denoted by $\nabla^2 f$, is a symmetric (0,2)-tensor defined by

\[\nabla^2_{X, Y} f \equiv XYf-(\nabla_X Y)f.\]

For a $(k, l)$ tensor $T$, we define a $(k, l+2)$ tensor $\nabla^2T$ by

\[\nabla^2_{X, Y} T \equiv \nabla_X\nabla_Y T-\nabla_{\nabla_X Y} T.\] This is no longer symmetric in $X$ and $Y$ and the commutator is the Riemannian curvature tensor.

New tensors can be constructed by contracting old ones. We can pick any orthonormal basis $\{e_A\}_{A=1}^N$ of the tangent bundle for contractions. The contracted tensor will be independent of the choice of such $e_A$. Notice that the summation convention, i.e. repeated indexes are summed, is followed throughout the article.
From the Riemannian curvature tensor, we define the Ricci curvature $Ric (X, Y)\equiv R(X, e_A, Y, e_A)$ and the scalar curvature $R\equiv Ric (e_A, e_A)$.

The Laplacian of $f$ is defined to be the trace of
$\nabla^2 f$, i.e.
\[\Delta f \equiv tr (\nabla^2 f)=\nabla^2_{e_A, e_A} f=\langle \nabla_{e_A} \nabla f, e_A\rangle
=(\nabla_{e_A} df)(e_A).\]

The rough Laplacian of $T$ is then
\[\Delta T \equiv tr \nabla^2 T=\nabla^2_{e_A, e_A} T.\]

We can always use  normal coordinates to simplify the
expression of a tensor at a given point $p$. Given any orthonormal basis of $T_p M$, there exists
a coordinate $\{y^A\}$ near $p$ and such that $e_A=\frac{\partial}{\partial y^A}$  and $\nabla_{e_A} e_B=0$ at the point $p$.

\subsection{Immersed submanifolds and the second fundamental forms}

Given an immersion $F:\Sigma\rightarrow M$ of an $n$-dimensional smooth manifold $\Sigma$ into $M$. Suppose  $\{y^A\}_{A=1}^N$ is a local coordinate system on $M$ with metric tensor $\Lambda_{AB}=\langle \frac{\partial}{\partial y^A}, \frac{\partial}{\partial y^B}\rangle$.  Because $F$ is an immersion, the tangent space of $\Sigma $ at $p$, $T_p\Sigma$ can be identified with $F_*T_p\Sigma$, the vector subspace of $ T_{F(p)}M$ spanned by
$\{ F_*(\frac{\partial}{\partial x^i})\}_{i=1}^n$.  We denote $F^A=y^A\circ F$, then $F_*(\frac{\partial}{\partial x^i})=\frac{\partial F^A}{\partial x^i}\frac{\partial}{\partial y^A}$. The metric $\Lambda_{AB}$ induces a Riemannian metric $g_{ij}$ on $\Sigma$ defined by
\[ g_{ij}=\langle F_*(\frac{\partial}{\partial x^i}), F_*(\frac{\partial }{\partial x^j})\rangle
=\frac{\partial F^A}{\partial x^i}\frac{\partial F^B}{\partial x^j} \Lambda_{AB}.\]

The orthogonal complement of $F_*T_p\Sigma$ in $T_{F(p)} M $, denoted by $N_p \Sigma$ is called the normal space of $\Sigma$ in $M$ at
the point $p$. Given a vector $V\in T_{F(p)}M$, we can define the projections of $V$ onto $F_*T_p\Sigma$ and $N_p\Sigma$, called
the tangential component $V^\top$ and the normal component $V^\perp$, respectively. They are
\[ V^\top=\langle V, F_*(\frac{\partial }{\partial x^i}) \rangle g^{ij} F_*(\frac{\partial}{\partial x^j})\in F_* T_p\Sigma, \] and
\[V^\perp=V-V^\top\in N_p \Sigma. \]

The Levi-Civita connection
$\nabla^M$ on $M$ induces a connection on $\Sigma$ by
\[\nabla^\Sigma_X Y\equiv (\nabla^M_X Y)^\top,\] for $X$ and $Y$ tangent to $\Sigma$.  It is not hard to check that $\nabla^\Sigma$ is the Levi-Civita connection of the induced metric $g_{ij}$ on $\Sigma$.

Given  a normal vector field $V$ along $\Sigma$.
The induced connection on the normal bundle, the normal connection, is defined to be

\[\nabla^\perp_X V \equiv (\nabla^M_X V)^\perp.\]

The curvature of the normal bundle is
\[R^\perp(X, Y) V \equiv -\nabla^\perp_X\nabla^\perp_Y V+\nabla^\perp_Y\nabla^\perp_X V+
\nabla^\perp_{[X, Y]} V.\]

The second fundamental form is defined to be

\[\mbox{II}(X, Y)\equiv (\nabla^M_X Y)^\perp,\] as a section of the tensor bundle $T^*\Sigma\otimes
T^*\Sigma\otimes N\Sigma$. Choose  orthonormal bases $\{e_i\}_{i=1}^n$ for $T_p\Sigma$
and $\{e_\alpha\}_{\alpha=n+1}^N$ for $N_p\Sigma$, the components of the second fundamental form are

\[h_{\alpha ij}=\langle\nabla^M_{e_i} e_j, e_\alpha\rangle
=\langle \mbox{II}(e_i, e_j), e_\alpha\rangle.\]

The mean curvature vector is the trace of the second fundamental form,

\[{H}\equiv \mbox{II}(e_i, e_i)=g^{ij}(\nabla^M_{F_*(\frac{\partial}{\partial x^i})}F_*(\frac{\partial}{\partial x^j}))^\perp.\]

Suppose the ambient space $M$ is $\R^N$ and let \[F=(F^1(x^1,\cdots, x^n), \cdots, F^N(x^1,\cdots, x^n))\]
be the position vector of the immersion. In this case,
we have
\[\Delta^{\Sigma} F=\Delta^{\Sigma}(F^1, \cdots, F^N)=(\Delta^\Sigma F^1,\cdots, \Delta^{\Sigma}F^N)={H},\] where $\Delta^\Sigma$ is the Laplace operator with respect to the induced metric $g_{ij}$ on $\Sigma$.

Denote the Riemannian curvature tensor of $g_{ij}$ by $R^\Sigma$. We recall the Gauss equation
\begin{equation}\label{gauss}\langle R^\Sigma (X, Y) Z, W\rangle=\langle R^M(X,Y) Z, W\rangle+\langle \mbox{II}(X, Z),\mbox{II} (Y, W)\rangle
-\langle \mbox{II}(X, W), \mbox{II} (Y, Z)\rangle,\end{equation} and the Codazzi equation
\begin{equation}\label{coda}(\nabla_X \mbox{II})(Y, Z)-(\nabla_Y \mbox{II})(X, Z)=-(R^M(X, Y)Z)^\perp.\end{equation} Here $\nabla$ on $\mbox{II}$ is the connection on the tensor bundle $T^*\Sigma\otimes
T^*\Sigma\otimes N\Sigma$, thus \begin{equation}\label{conn_second}(\nabla_X \mbox{II})(Y, Z)=\nabla^\perp_X \mbox{II}(Y, Z)-\mbox{II}(\nabla^\Sigma_X Y, Z)
-\mbox{II}(Y, \nabla^\Sigma_X, Z).\end{equation}

\subsection{First variation formula}
Recall the volume of $\Sigma$ with respect to the metric $g$ is

\[Vol(\Sigma)=\int_{\Sigma} \sqrt{\det g_{ij}}dx^1\wedge\cdots\wedge dx^n\]
Take a normal vector field $V$ along $\Sigma$ and consider a family of immersion $F:\Sigma\times [0, \epsilon)\rightarrow M$ such that $\frac{\partial F}{\partial s}|_{s=0} =V$. Denote $\Sigma_s =F(\Sigma, s)$.

To consider the change of volume
\[\frac{d}{ds}|_{s=0} Vol(\Sigma_s)=\frac{d}{ds}|_{s=0}\int_\Sigma\sqrt{\det g_{ij}(s)}dx^1\wedge\cdots \wedge dx^n,\]
we compute
\[\frac{d}{ds}\sqrt{\det g_{ij}}=\frac{1}{2}(\det g_{ij})^{-\frac{1}{2}}(\frac{d}{ds} g_{ij})(\det g_{ij}),\] and
\[\frac{d}{ds} g_{ij}=-2\langle V, \nabla^M_{F_*(\frac{\partial }{\partial x^i})}, F_*(\frac{\partial}{\partial x^j})\rangle.\]

Therefore,
\[\frac{d}{ds}|_{s=0} Vol (\Sigma_s)=-\int_{\Sigma} \langle V, {H}\rangle \sqrt{\det g_{ij}}dx^1\wedge\cdots \wedge dx^n\] and the mean curvature vector is formally the gradient of the negative of the volume functional on the space of submanifolds.

\section{Mean curvature flow}
\subsection{The equation}

Let $F_0: \Sigma \rightarrow M$ be a smooth immersion.
Consider $F: \Sigma \times [0, T) \rightarrow M$ that satisfies
\[ \left\{ \begin{array}{l}
\frac{\partial F}{\partial t} = H \smallskip \\
F(\cdot,0) = F_0(\cdot)
\end{array} \right. \]
$F$ is called the mean curvature flow of $F_0(\Sigma)$ (in the normal direction). By the first variation formula, we have
\[\frac{d}{dt}\sqrt{\det g_{ij}}=-|{H}|^2\sqrt{\det g_{ij}}.\] Denote $F(\cdot, t) $ by $F_t:\Sigma\rightarrow M$. Because  $\det g_{ij}>0$ at $t=0$, $F_t$ remains an immersion as long as $|{H}|^2$ is bounded.

It is not hard to established the short time existence of the flow in the case of $\Sigma$ is compact.
Notice that the mean curvature flow is not strictly parabolic, suppose $\phi_t$ , $t\in [0, T)$
is a family of diffeomorphism of $\Sigma$. Consider $\bar{F}_t=F_t\circ \phi_t: \Sigma
\rightarrow M$, the image of $\bar{F}_t$ is the same as that of $F_t$. But $\frac{\partial \bar{F}_t}{\partial t}$ has both the tangential and normal components. Therefore
the general form of mean curvature flow is
\[(\frac{\partial F}{\partial t})^\perp={H}.\]

On the other hand, given a general mean curvature flow $(\frac{\partial \bar{F}}{\partial t})^\perp={H}$,
we can always find a reparametrization $\phi_t$ so that $F_t=\overline{F}_t\circ \phi_t$ satisfies
$\frac{\partial F}{\partial t}={H}$.

Suppose $M=\R^{n+1}$ and $\Sigma$ is the
the graph of a function $f(x^1,\cdots, x^n)$ on $\{x^{n+1}=0\}$. We can represent the hypersurface
in the parametric form $F(x^1, \cdots, x^n, t)=(x^1, \cdots, x^n, f(x^1, \cdots, x^n, t))$. Denote by $\nu$ the
upward unit normal $\frac{(-\nabla f, 1)}{\sqrt{1+|\nabla f|^2}}$ on the graph, by the mean curvature flow equation $(\frac{\partial F}{\partial t})^\perp={H}$, we
derive \[\langle \frac{\partial F}{\partial t}, \nu\rangle\nu =H. \] The mean curvature flow of the graph of $f$ becomes a quasi-linear parabolic equation for $f$:
\[\frac{1}{\sqrt{1+|\nabla f|^2}}\frac{\partial f}{\partial t}=div(\frac{\nabla f}{\sqrt{1+|\nabla f|^2}}).\] This equation was studied extensively in \cite{eh1} and \cite{eh2}.

In general, suppose  $F:\Sigma\times [0, T)\rightarrow \R^N$ and $F(x^1,\cdots, x^n, t)$ is
a mean curvature flow. As ${H}=(g^{ij} \frac{\partial^2 F}{\partial x^i\partial x^j})^\perp$,
the mean curvature flow equation is
\[\frac{\partial F}{\partial t}=(g^{ij} \frac{\partial^2 F}{\partial x^i \partial x^j})^\perp.\]

Given any vector $V\in \R^N$, we compute
\[V^\perp=V-V^\top=V-\langle V, \frac{\partial F}{\partial x^k}\rangle g^{kl}\frac{\partial F}{\partial x^l}=V^A\frac{\partial}{\partial y^A}-V^B\frac{\partial F^B}{\partial x^k}g^{kl} \frac{\partial F^C}{\partial x^l}\frac{\partial}{\partial y^C}\]

Therefore, the mean curvature flow can be expressed as
\[\frac{\partial F^B}{\partial t}=g^{ij} \frac{\partial^2 F^A}{\partial x^i\partial x^j}
(\delta_{AB}-\frac{\partial F^B}{\partial x^k} g^{kl} \frac{\partial F^A}{\partial x^l}).\]
\subsection{Finite time singularity}
Suppose ${F}_0$ is the standard imbedding of the unit sphere $S^n$ in $\R^N$ and  $F(x, t)=\phi(t){F}_0(x)$ is a mean curvature flow in the normal direction. As
\[\frac{\partial F}{\partial t}=\phi'(t) {F}_0(x)\]
and the mean curvature vector of $F(x, t)$ is
\[{H} (x, t)=\frac{n}{\phi(t)} (-{F}_0(x)),\] the mean curvature flow is reduced to the ordinary differential equation \[\phi'(t)=-\frac{n}{\phi(t)}.\]  We solve $\phi(t)=\sqrt{r_0^2-2nt}$ where $r_0$ is the radius of $F(x, 0)$. Therefore,
\[F(x,t) =\sqrt{r_0^2-2nt} F_0(x).\] $F(x,t)$ shrinks to a point at $t=\frac{r_0^2}{2n}$. This is the limiting behavior of the mean curvature flow of any compact convex hypersurface in $\R^N$, see \cite{hu1} and \cite{hu2}.

Indeed, it is very easy to show the mean curvature flow of any compact submanifolds in $\R^N$ must develop finite time singularity. Let $F:\Sigma \times[0, T)\rightarrow \R^N$ be such a flow for a compact $n$-dimensional $\Sigma$. By the formula of the mean curvature vector, we have $\frac{d}{dt} F=\Delta_t F$, where $\Delta_t$ is the Laplace operator with respect to the induced metric on $F_t(\Sigma)$. We compute

\begin{equation}\frac{d}{dt}|F|^2=\Delta_t |F|^2-2n.\end{equation} Therefore $|F|^2+2nt$ satisfies the heat equation:

\begin{equation}\frac{d}{dt}(|F|^2+2nt)=\Delta_t (|F|^2+2nt).\end{equation} Suppose $F_0(\Sigma)$ is contained in a sphere $S^{N-1}(r_0)$ of radius $r_0$ centered at the origin and thus $|F|^2\leq r_0^2$ at $t=0$. It follows from the maximum principle that $|F|^2+2nt\leq r_0^2$ must be preserved as long as the flow exists smoothly. In particular, the flow must develop singularity before the time $t=\frac{r_0^2}{2n}$.

\section{Blow-up analysis}
\subsection{Backward heat kernel and monotonicity formula}

Fix $y_0 \in \mathbb{R}^{N}$ and $t_0\in \R$ and consider the backward heat kernel at $(y_0, t_0)$,
\begin{equation} \label{backward}\rho_{y_0,t_0}(y,t) = \frac{1}{(4 \pi (t_0-t))^{n/2}} \exp \left( -\frac{|y-y_0|^2}{4(t_0-t)} \right)\end{equation}
defined on $\mathbb{R}^{N} \times (-\infty,t_0)$.

Suppose $F: \Sigma \times [0,t_0) \rightarrow \R^N$ is a mean curvature flow, then Huisken's monotonicity formula \cite{hu2} says
\begin{equation}\label{mono_eucl} \frac{d}{dt} \int_{\Sigma_t} \rho_{y_0,t_0} d\mu_t =
-\int_{\Sigma_t} \rho_{y_0,t_0} \left| \frac{F^\perp}{2(t_0-t)}+H \right| ^2 d\mu_t \leq 0.\end{equation}
Hence $\lim_{t \rightarrow t_0} \int_{\Sigma_t} \rho_{y_0,t_0} d\mu_t$ exists.

 For a general ambient manifold $M$ of dimensional $N$, we fix an isometric embedding $ i:M \rightarrow \mathbb{R}^{N+m}$. Given an immersion $F: \Sigma \rightarrow M$, we consider $\bar{F}=i\circ F$ as an immersed submanifold in $\R^{N+m}$. Denote by $H$ the mean curvature vector of $\Sigma$ with respect to $M$, i.e. \[H=(\nabla^M_{e_i} e_i)^\perp
\in TM\] for an orthonormal basis $\{e_i\}$ of $T\Sigma$. Denote by $\bar{H}$ the mean curvature vector of $\Sigma$ with respect to $\R^{N+m}$, i.e. \[\bar{H}=(\nabla^{\R^{N+m}}_{e_i} e_i)^\perp\in T\R^{N+m}.\] Therefore,
\begin{equation}\label{eq_E}\bar{H}-H=\mbox{II}^M(e_i, e_i)=-E,\end{equation} where $\mbox{II}^M$ is the second fundamental form of $M$ in $\R^{N+m}$. $|E|$ is bounded if $M$ has bounded geometry. The equation of the mean curvature flow becomes
\begin{equation}\label{mcf}\frac{\partial \bar{F}}{\partial t}=\bar{H}+E.\end{equation}
The following monotonicity formula for a
general ambient manifold  was derived by White in \cite{wh2}.

\begin{prop}
For the mean curvature flow (\ref{mcf}) on $[0, t_0)$, if we isometrically embed $M$ into $\R^{N+m}$ and take the backward heat kernel $\rho_{y_0, t_0}$ for $y_0\in \R^{N+m}$, then
\begin{equation}\label{mono_gen} \frac{d}{dt} \int_{\Sigma_t} \rho_{y_0,t_0} d\mu_t \leq
      C- \int_{\Sigma_t} \left| H+\frac{1}{2(t_0-t)}\bar{F}^{\perp}+\frac{E}{2} \right| ^2 \rho_{y_0,t_0} d\mu_t \end{equation} for some constant $C$. Here $\bar{F}^\perp$ is the component of the position vector $\bar{F}\in \R^{N+m}$ in the normal space of $\Sigma_t$ in $\R^{N+m}$.
\end{prop}

\begin{defn}
The density at $(y_0, t_0)$ is defined to be $\Theta(y_0,t_0) = \lim_{t \rightarrow t_0} \int_{\Sigma_t} \rho_{y_0,t_0} d\mu_t $.
\end{defn}

We will use the following pointwise formula which can be found in \cite{wa1}
\begin{equation}\label{bh}
\begin{split}
\frac{d}{dt}\rho_{y_0,t_0}&=-
\Delta^{\Sigma_t}\rho_{y_0,t_0}
-\rho_{y_0,t_0}( \frac{|\bar{F}^\perp|^2}{4(t_0-t)^2}
+\frac{\bar{F}^\perp \cdot \bar{H}}{t_0-t}
+\frac{\bar{F}^\perp \cdot E}{2(t_0-t)})
\end{split}
\end{equation}

The minus sign in front of the Laplacian in equation (\ref{bh})
indicates that $\rho_{y_0, t_0}$ satisfies the backward heat
equation.

\subsection{Synopsis of singularities}
\begin{defn}Given a mean curvature flow $F:\Sigma \times [0, t_0)\rightarrow M$, the image $\mathcal{F}$  of the map
$$ \begin{array}{cccc}
F \times 1: &\Sigma \times [0,t_0) &\rightarrow &M \times \mathbb{R} \\
&(x,t) &\mapsto &(F(x,t),t)
\end{array}$$ in $M\times \R$ is called the space-time track of the flow. If $M$ is isometrically embedded into $\R^{N+m}$, we can take the image in $\R^{N+m}\times \R$ as the space-time track.
\end{defn}

\begin{defn}
The parabolic dilation of scale $\lambda>0$ at $(y_0, t_0)$ is given by
$$ \begin{array}{lccl}
D_\lambda: &\mathbb{R}^{N+m} \times \mathbb{R} &\rightarrow &\mathbb{R}^{N+m} \times \mathbb{R} \\
& (y,t) &\mapsto &(\lambda (y-y_0), \lambda^2( t-t_0)) \smallskip \\
& t \in [0,t_0) &\mapsto & s = \lambda^2(t-t_0) \in [-\lambda^2 t_0,0)
\end{array} $$
\end{defn}
Set $\Sigma_s^\lambda=D_\lambda (\Sigma_t-y_0)=\lambda (\Sigma_{t_0+\frac{s}{\lambda^2} }-y_0)$ for $t=t_0+\frac{s}{\lambda^2}$. If $M = \mathbb{R}^N$, the mean curvature flow equation is preserved by $D_\lambda$. That is, if $\mathcal{F}$ is the spacetime track of a mean curvature flow, so is $D_\lambda \mathcal{F}$.

\vskip 10pt
\includegraphics[height=4cm]{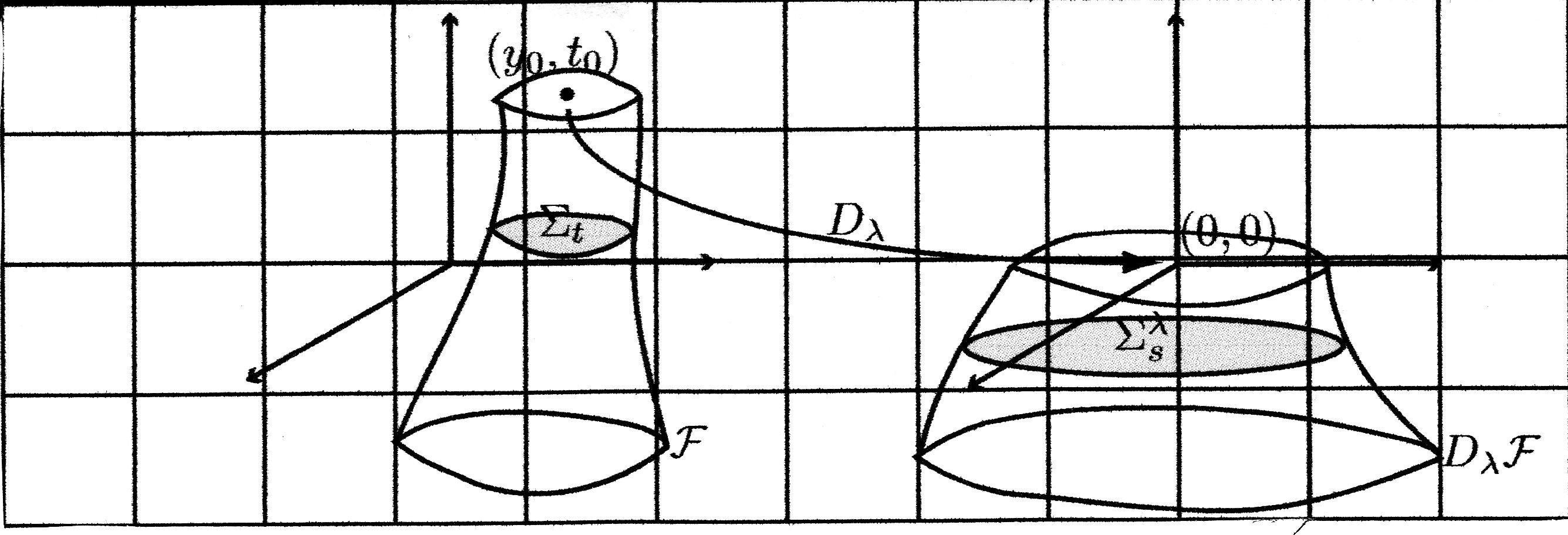}

\vskip 10pt

\begin{prop}
$\int_{\Sigma_t} \rho_{y_0,t_0} d \mu_t$ is invariant under the parabolic dilation. That is to say,
$$ \int_{\Sigma_t} \rho_{y_0,t_0} d \mu_t = \int_{\Sigma^\lambda_s} \rho_{0,0} d \mu^\lambda_s $$
where $\mu^\lambda_s$ is the volume form of the $\Sigma_s^\lambda$.
\end{prop}

The mean curvature flow extends smoothly to $t_0 + \delta$ if the second fundamental form is bounded
$\sup_{\Sigma_t} |\mbox{II}| < C$ as $t\rightarrow t_0$.
Therefore if a singularity is forming at $t_0$, then $\sup_{\Sigma_t} |\mbox{II}| \rightarrow \infty$ as $t\rightarrow t_0$.
\begin{defn}[type I singularity]
A singularity $t=t_0$ is called a type I singularity if there is a $C>0$ such that
\[ \sup_{\Sigma_t} |\text{\mbox{II}}|^2 \leq \frac{C}{t_0-t}\] for all $t<t_0$.
\end{defn}

The shrinking sphere $F(x,t)=\sqrt{r_0^2-2t} \,F_0(x)$ is a type I singularity at $t_0=\frac{r_0^2}{2n}$.
 For any $t<t_0$, the radius of the sphere is $r= \sqrt{r_0^2-2nt}$ and thus the principal curvature $k$ satisfies
$$ k^2 = \frac{1}{r^2} = \frac{1}{r_0^2-2nt} = \frac{1}{2n(t_0-t)}. $$

Consider the parabolic dilation of the spacetime  track of a type I singularity, we calculate
$$ |\mbox{II}|^2(\Sigma_s^\lambda) = \frac{1}{\lambda^2}|\mbox{II}|^2(\Sigma_t)
     = \frac{-1}{s}(t_0-t)|\mbox{II}|^2(\Sigma_t) \leq \frac{-1}{s} C \mbox{ , for  } s \in [-\lambda^2 t_0, 0).$$
Hence for any fixed $s$, the second fundamental form $|\mbox{II}|^2(\Sigma_s^\lambda)$ is bounded.
By Arzel\`a-Ascoli theorem, on any compact subset of space time,
there exists smoothly convergent subsequence of $\{ D_\lambda \mathcal{F} \}$ as $\lambda\rightarrow \infty$.
We thus have:
\begin{prop}
If there is  a type I singularity at $t=t_0$, there exists a subsequence $\{ \lambda_i \}$ such that
$D_{\lambda_i} \mathcal{F} \rightarrow \mathcal{F}_\infty$, the spacetime track of a smooth flow that exists on $(-\infty, 0)$.
\end{prop}

After the parabolic dilation,  with $t= t_0+\frac{s}{\lambda^2}$, the monotonicity formula (\ref{mono_gen}) becomes
\[\frac{d}{ds} \int_{\Sigma_s^\lambda} \rho_{0,0} d\mu_s^\lambda \leq \frac{C}{\lambda^2} - \int_{\Sigma_s^\lambda} \left| H_s^\lambda+\frac{-1}{2s}(\bar{F}_s^\lambda)^{\perp}
       +\frac{E}{2\lambda} \right|^2 \rho_{0,0} d\mu_s^\lambda.\]
Consider the $s_0<0$ slice and integrate both sides from $s_0-\tau$ to $s_0$ for $\tau>0$ and $\lambda$ large:
\begin{align*}
&\int_{s_0-\tau}^{s_0} \int_{\Sigma_s^{\lambda}} \left| H_s^{\lambda}-\frac{1}{2s}(\bar{F}_s^{\lambda})^{\perp}
       +\frac{E}{2\lambda} \right|^2 \rho_{0,0} \, d\mu_s^{\lambda} ds \\
\leq &\frac{C \tau}{\lambda^2} - \int_{\Sigma_{s_0}^{\lambda}} \rho_{0,0} \, d\mu_{s_0}^{\lambda}
       + \int_{\Sigma_{s_0-\tau}^{\lambda}} \rho_{0,0} \, d\mu_{s_0-\tau}^{\lambda}
\end{align*}
Take $\lambda \rightarrow \infty$. Because
$ \int_{\Sigma_{s_0}^{\lambda}} \rho_{0,0} \, d\mu_{s_0}^{\lambda}
=\int_{\Sigma_t} \rho_{y_0,t_0} \, d\mu_t \mbox{ where } t = t_0 + \frac{s_0}{\lambda^2} $,
$$\lim_{\lambda \rightarrow \infty} \int_{\Sigma_{s_0}^{\lambda}} \rho_{0,0} \, d\mu_{s_0}^{\lambda}
=\lim_{t \rightarrow t_0} \int_{\Sigma_t} \rho_{y_0,t_0} \, d\mu_t=\lim_{\lambda \rightarrow \infty} \int_{\Sigma_{s_0-\tau}^{\lambda}} \rho_{0,0} \, d\mu_{s_0-\tau}^{\lambda}.$$
Therefore,
$$ \lim_{\lambda \rightarrow \infty} \int_{s_0-\tau}^{s_0} \int_{\Sigma_s^{\lambda}} \left| H_s^{\lambda}-
      \frac{1}{2s}(\bar{F}_s^{\lambda})^{\perp}+\frac{E}{2\lambda} \right|^2 \rho_{0,0} \, d\mu_s^{\lambda} ds = 0$$
Therefore if $t_0$ is a type I singularity, a subsequence $D_{\lambda_i} \mathcal{F} \rightarrow \mathcal{F}_\infty$ smoothly and $H = \frac{1}{2s}F^\perp$ on $\mathcal{F}_\infty$ for $s \in (s_0-\tau, s_0)$. Take $s_0\rightarrow 0$ and $\tau\rightarrow \infty$, we obtain Huisken's theorem \cite{hu2}:
\begin{thm}
If the singularity is of type I, then there exists a subsequence $\lambda_i$ such that $D_{\lambda_i} \mathcal{F} \rightarrow \mathcal{F}_\infty $ smoothly and $H = \frac{1}{2s}F^\perp, -\infty<s<0$ on $\mathcal{F}_\infty$.
\end{thm}

Let $F$ be a solution of the mean curvature flow  that satisfies $H(x,s) = \frac{1}{2s}F^\perp(x,s)$ on $-\infty<s<0$ then $F$ must be of the form $F(x,s) = \sqrt{-s}F(x)$. Even if the singularity is not of type I (so the scaled $|\mbox{II}|$ is not necessarily bounded),
we still get weak solution with $H = \frac{1}{2s} F^\perp$ (see \cite{il2}). A theorem of Huisken \cite{hu2} shows any compact smooth convex hypersurface in $\mathbb{R}^N$ satisfying  $H = - F^\perp$ is the unit sphere.

\section{Applications to deformations of symplectomorphisms of Riemann surfaces}

\subsection{Introduction}
We apply the mean curvature flow to obtain a result of the deformation retract of symplectomorphism groups of Riemann surfaces.
\begin{thm} \label{mainthm}
Let $\Sigma^1$ and $\Sigma^2$ be two compact closed Riemann surfaces with metrics of the same constant
curvature $c$. Let $\omega_1$  and $\omega_2$ be the volume forms of $\Sigma^1$ and
$\Sigma^2$, respectively.
Consider a map $f: \Sigma^1 \rightarrow \Sigma^2$ that satisfies $f^* \omega_2 = \omega_1$, i.e.  $f$ is an area-preserving map or a symplectomorphism.
Denote by $\Sigma_t$ the mean curvature flow of the graph of $f$ in $M = \Sigma^1 \times \Sigma^2$. We have
\begin{enumerate}
\item $\Sigma_t$ exists smoothly for all $t>0$ and converges smoothly to $\Sigma_\infty$ as $t \rightarrow \infty$.
\item Each $\Sigma_t$ is the graph of a symplectomorphism $f_t:\Sigma^1\rightarrow \Sigma^2$ and $f_t$ converges smoothly to a symplectomorphism $f_\infty:\Sigma^1\rightarrow \Sigma^2$ as $t\rightarrow \infty$.
\end{enumerate}
Moreover,
$$ f_\infty\,\, \text{is}\,\, \left\{ \begin{array}{ll}
\mbox{an isometry} & \mbox{if } c >0 \\
\mbox{a linear map} & \mbox{if } c = 0 \\
\mbox{a harmonic diffeormophism} & \mbox{if } c<0
\end{array} \right. $$
\end{thm}
The existence of such minimal symplectic maps in the $c<0$ case was proved by Schoen in \cite{sc} (see also \cite{le}). This theorem was proved in \cite{wa2} and \cite{wa3}. A different proof in the case $c\leq 0$ was proved by Smoczyk \cite{sm2} assuming an extra angle condition.

Both $\omega_1$ and $\omega_2$ can be extended to parallel forms on $M$. Denote $\omega_1-\omega_2$ by $\omega'$ and $\omega_1+\omega_2$ by $\omega''$, which are parallel forms on $M$ as well. That $\Sigma_t$ is the graph of a symplectomorphism $f_t$ can be characterized by $*_{\Sigma_t}\omega'=0$, $*_{\Sigma_t}\omega''>0$ where $*_{\Sigma_t}$ is the Hodge star operator on $\Sigma_t$. To prove that $f_t$ being a symplectomorphism is preserved, it suffices to show both conditions are preserved along the mean curvature flow. In the following, we compute the evolution equations of these quantities.

\subsection{Derivation of evolution equations}

In this section, we assume $M$ is a four-dimensional Riemannian manifold and $\Sigma$ is a closed Riemann surface embedded in $M$.
\begin{lem}
Let $\omega$ be a parallel 2-form on $M$
and $\eta=*_{\Sigma}\omega$, then
\begin{align*}
\Delta^{\Sigma}\eta &=-\eta|\mbox{II}|^2+\omega(\nabla^\perp_{e_1}H,e_2)+\omega(e_1,\nabla^\perp_{e_2}H) \\
&-\omega((R^M(e_1,e_k)e_k)^\perp,e_2)-\omega(e_1,(R^M(e_2,e_k)e_k)^\perp)
	+2\omega(\mbox{II}(e_k,e_1),\mbox{II}(e_k,e_2))
\end{align*}
where $\{e_i\}=\{e_1,e_2\}\in T_p\Sigma$ and $\{e_\alpha\}=\{e_3,e_4\}\in N_p\Sigma$ are oriented orthonormal bases and $\mbox{II}(e_i,e_j)=(\nabla^M_{e_i}e_j)^\perp$ is the second fundamental form of $\Sigma$.
\end{lem}

\begin{proof} The right hand side is independent of the choice of frames and thus a tensor (it does depend on the orientation though).
At any $p\in\Sigma$, we can choose an orthonormal frame $\{e_1,e_2\}$ such that
$\nabla^{\Sigma}_{e_i}e_j =0$ at $p$. It is not hard to check that $\Delta^\Sigma \eta=*\Delta^\Sigma \omega$ where the $\Delta^\Sigma$ on the right hand side is the rough-Laplacian on two forms. Therefore $\Delta^\Sigma \eta$ is equal to
\begin{align*}
(\nabla^{\Sigma}_{e_k}\nabla^{\Sigma}_{e_k}\omega)(e_1, e_2)
=&e_k\big( (\nabla^\Sigma_{e_k}\omega)(e_1,e_2) \big)-(\nabla^\Sigma_{e_k}\omega)(\nabla^\Sigma_{e_k}e_1,e_2)
	-(\nabla^\Sigma_{e_k}\omega)(e_1,\nabla^\Sigma_{e_k}e_2) \\
=&e_k\big( e_k(\omega(e_1,e_2))-\omega(\nabla^\Sigma_{e_k}e_1,e_2)-\omega(e_1,\nabla^\Sigma_{e_k}e_2) \big).
\end{align*}

Since $\omega$ is parallel on $M$, we derive
\begin{align*}
\Delta^\Sigma \eta= &e_k\big(\omega(\nabla^M_{e_k}e_1,e_2)+\omega(e_1,\nabla^M_{e_k}e_2)
	-\omega(\nabla^\Sigma_{e_k}e_1,e_2)-\omega(e_1,\nabla^\Sigma_{e_k}e_2) \big) \\
=&e_k\big( \omega(\mbox{II}(e_k,e_1),e_2)+\omega(e_1,\mbox{II}(e_k,e_2)) \big). \\
\end{align*}
Applying $\nabla^M\omega=0$ again, the last expression can be decomposed into
\begin{align*}
&\omega(\nabla^M_{e_k}\mbox{II}(e_k,e_1),e_2)+\omega(\mbox{II}(e_k,e_1),\nabla^M_{e_k}e_2)+\omega(\nabla^M_{e_k}e_1,\mbox{II}(e_k,e_2))+\omega(e_1,\nabla^M_{e_k}\mbox{II}(e_k,e_2)) \\
=&\omega((\nabla^M_{e_k}\mbox{II}(e_k,e_1))^\top,e_2)+\omega(\nabla^\perp_{e_k}\mbox{II}(e_k,e_1),e_2) \\
	&+\omega(e_1,(\nabla^M_{e_k}\mbox{II}(e_k,e_2))^\top)+\omega(e_1,\nabla^\perp_{e_k}\mbox{II}(e_k,e_2))
	+2\omega(\mbox{II}(e_k,e_1),\mbox{II}(e_k,e_2)).
\end{align*}
We compute
\[\omega((\nabla^M_{e_k}\mbox{II}(e_k,e_1))^\top,e_2) =\langle \nabla^M_{e_k}\mbox{II}(e_k,e_1), e_1\rangle \eta=-\sum_{k=1}^2|\mbox{II}(e_k,e_1)|^2\eta\] and likewise for the other term $\omega(e_1,(\nabla^M_{e_k}\mbox{II}(e_k,e_2))^\top)$.
On the other hand,
\[\omega(\nabla^\perp_{e_k}\mbox{II}(e_k,e_1),e_2)=\omega((\nabla_{e_k}\mbox{II})(e_k,e_1),e_2)=\omega((\nabla_{e_k}\mbox{II})(e_1,e_k),e_2)\] where the formula for the connection on $\mbox{II}$ (\ref{conn_second}) is used.
Applying the Codazzi equation (\ref{coda}), we see this is equal to
\[\omega((\nabla_{e_1}\mbox{II})(e_k,e_k),e_2)-\omega((R^M(e_k,e_1)e_k)^\perp,e_2)=\omega(\nabla^\perp_{e_1} H,e_2)-\omega((R^M(e_k,e_1)e_k)^\perp,e_2).\] The formula is proved.
	
\end{proof}

We compute  the evolution equation of $*_{\Sigma_t}\omega$ along a mean curvature flow.
\begin{prop}
Let $\Sigma_t$ be mean curvature flow of a closed embedded Riemann surface in $M$. Let $\omega$ be a parallel 2-form on $M$
and $\eta=*_{\Sigma_t}\omega$, then
\begin{align*}
(\frac{d}{dt}-\Delta^{\Sigma_t})\eta &= |\mbox{II}|^2\eta-2\omega(\mbox{II}(e_k,e_1),\mbox{II}(e_k,e_2)) \\
&\quad +\omega((R^M(e_1,e_k)e_k)^\perp,e_2)+\omega(e_1,(R^M(e_2,e_k)e_k)^\perp)
\end{align*}
\end{prop}
\begin{proof}  We parametrize $\Sigma_t$ by $F:\Sigma\times[0, T)\rightarrow M$ with $\frac{\partial F}{\partial t}=H$ and fix a local oriented coordinate system $x^1, x^2$ on $\Sigma$.
Write $\eta=\frac{1}{\sqrt{\det g_{ij}}}\omega(F_*(\frac{\partial }{\partial x^1}), F_*(\frac{\partial }{\partial x^2}))$ and compute
$$ \frac{d}{dt}\eta=\omega((\nabla_{e_1}H)^\perp,e_2)-\omega((\nabla_{e_2}H)^\perp,e_1). $$
Combine this with the formula of $\Delta^{\Sigma_t}\eta$ in Lemma 4.1, we obtain the desired formula.
\end{proof}

\begin{defn}
Given a Riemannian manifold $M$, we say a two-form $\omega$ is of {\textit K\"ahler type} on $M$ if  $\omega(X, Y)=\langle JX, Y\rangle$ for a parallel complex structure $J (J^2=-I)$ that is compatible with the metric, i.e
$\langle JX, JY\rangle=\langle X, Y\rangle$.
\end{defn}
Both $\omega'$ and $\omega''$ in the previous section are of K\"ahler type on $M=\Sigma^1\times\Sigma^2$.

\begin{lem}\label{Ricci}
If $\omega$ is of K\"ahler type on $M$ and $J$ is the corresponding complex structure, then
\[\omega((R^M(e_1,e_k)e_k)^\perp,e_2)+\omega(e_1, (R^M(e_2,e_k)e_k)^\perp)=(1-\eta^2)\mbox{Ric}^M(Je_1,e_2).\]\end{lem}
\begin{proof}
Recall the following identity
\[\mbox{Ric}^M(JX,Y)=\frac{1}{2}\sum_{A=1}^4 R^M(X,Y,e_A,Je_A)\] for any orthonormal frame $\{e_A\}_{A=1}^4$ of $T_pM$.
Denote $\omega(e_A, e_B)=\omega_{AB}$ and $R^M(e_A, e_B, e_C, e_D)=R_{ABCD}$.
When $\eta\neq\pm1$, we can choose an oriented orthonormal frame $\{e_1,e_2, e_3, e_4\}$ with $e_1, e_2\in T_p\Sigma$ and $e_3, e_4\in N_p\Sigma$, such that
\begin{equation}\label{frames}\omega_{12}=\omega_{34}=\eta, \omega_{13}=-\omega_{24}=\sqrt{1-\eta^2}, \text{and  } \omega_{23}=\omega_{14}=0.\end{equation}

We compute\[\mbox{Ric}^M(Je_1,e_2)=\frac{1}{2}R^M(e_1,e_2,e_A,Je_A)=\sum_{A<B} \omega_{AB} R_{12AB}.\]
Plug in (\ref{frames}), we obtain \[\mbox{Ric}^M(Je_1,e_2)=\eta(R_{1212}+R_{1234})+\sqrt{1-\eta^2}(R_{1213}-R_{1224}).\]
The $J$ invariance of the curvature operator implies $R^M(e_1,e_2,e_1,e_2)=R^M(e_1,e_2,Je_1,Je_2)$. Plug in (\ref{frames}) again, we obtain \[(1-\eta^2)(R_{1212}+R_{1234})=\eta\sqrt{1-\eta^2}(R_{1213}-R_{1224}).\] Therefore,
\[\mbox{Ric}^M(Je_1,e_2)=\frac{1}{\sqrt{1-\eta^2}}(R_{1213}-R_{1224}).\]
On the other hand,
\[\omega((R^M(e_1,e_k)e_k)^\perp,e_2)+\omega(e_1, (R^M(e_2,e_k)e_k)^\perp)
=R_{1kk4}\omega _{42}+R_{2kk3}\omega_{13}.\] The lemma is proved by recalling (\ref{frames}) again.
\end{proof}

The terms that involve the second fundamental form in Prop 4.1 can also be simplified using (\ref{frames})
\begin{align*}
-2\omega(h_{\alpha 1k}e_\alpha,h_{\beta 2k}e_\beta) &=-2h_{\alpha 1k}h_{\beta 2k}\omega(e_\alpha,e_\beta) \\
&=(-2h_{31k}h_{42k}+2h_{41k}h_{32k})\omega(e_3,e_4) \\
&=(-2h_{31k}h_{42k}+2h_{41k}h_{32k})\eta.\end{align*}

We have thus proved the following theorem.
\begin{thm}\label{eq_eta}
Suppose $\omega$ is of K\"ahler type on $M$ and $J$ is the corresponding complex structure, then
\begin{align*}
\big(\frac{d}{dt}-\Delta^{\Sigma_t} \big)\eta=((h_{31k}-h_{42k})^2+(h_{41k}+h_{32k})^2)\eta+(1-\eta^2)\mbox{Ric}^M(Je_1,e_2)
\end{align*}
In particular, if the Riemannian metric $g$ on $M$ is Einstein with $\mbox{Ric}^M=cg$,
\begin{equation}\label{eq_ke}\big( \frac{d}{dt}-\Delta^{\Sigma_t} \big)\eta
=((h_{31k}-h_{42k})^2+(h_{41k}+h_{32k})^2)\eta+c\eta(1-\eta^2).
\end{equation}
\end{thm}

\subsection{Long time existence}
We apply Theorem \ref{eq_eta} to $M=\Sigma^1\times \Sigma^2$ on which both $\omega'$ and $\omega''$ are of K\"ahler type.
The maximum principle implies $*_{\Sigma_t}\omega'=0$ is preserved along the mean curvature flow and thus $\Sigma_t$ remains Lagrangian. We remark that in \cite{sm2} Smoczyk proved that being Lagrangian is preserved by the mean curvature flow in any K\"ahler-Einstein manifold.

Set $\eta=*_{\Sigma_t}\omega''$ and apply the maximum principle, we see that $\eta>0$ is preserved, or $\Sigma_t$ remains the graph of a symplectomorphism.  Compare $\eta$ with the solution of the ordinary differential equation
\[ \frac{d}{dt}f=cf(1-f^2),\] we deduce
 \begin{equation}\label{eta_lower}\eta(x,t)\geq\frac{\alpha e^{ct}}{\sqrt{1+\alpha^2e^{2ct}}} \mbox{ , where } \alpha>0 \mbox{  is defined by} \frac{\alpha}{\sqrt{1+\alpha^2}}=\min_{\Sigma_0}\eta.\end{equation}

 Let $J'$ be the parallel complex structure associated with $\omega'$ and $\{ e_1, e_2 \}$ be an oriented orthonormal basis on $\Sigma_t$.
Take $e_3=J'(e_1), e_4=J'(e_2)$, then $\{ e_3, e_4 \}$ is an orthonormal basis for the normal space $N\Sigma_t$. The second fundamental has the following symmetries
\[\langle \mbox{II}(X, Y), J'Z\rangle=\langle \mbox{II}(X, J'Y), Z\rangle=\langle \mbox{II}(J'X, Y), Z\rangle.\] These symmetries imply
\begin{equation}\label{second_ineq}|H|^2\leq \frac{4}{3} |\mbox{II}|^2.\end{equation}
The evolution equation (\ref{eq_ke} ) becomes
\begin{equation}\label{eta_lg}\big( \frac{d}{dt}-\Delta^{\Sigma_t} \big)\eta
=\eta(2|\mbox{II}|^2-|H|^2)+c\eta(1-\eta^2).
\end{equation}
We fix an isometric embedding $i: M\rightarrow\mathbb{R}^N$.
Let $F$ be the mean curvature flow of $\Sigma$ in $M$ and denote $i\circ F$ by $\bar{F}$. We use the notation in \S 3.1.

For $y_0\in \R^N$, suppose $(y_0,t_0)$ is a singularity of the mean curvature. We consider the backward heat kernel $\rho_{y_0, t_0}$ defined in \S 3.1.

Denote $Q=1-\eta=$ and $P=-\eta(2|\mbox{II}|^2-|H|^2)-c\eta(1-\eta^2)$. Recalling equations (\ref{bh}) and (\ref{eta_lg}), we derive
\begin{align*}
\frac{d}{dt}\int_{\Sigma_t}Q\rho_{y_0,t_0}d\mu_t
&=\int_{\Sigma_t}(\Delta^{\Sigma_t} Q+P)\rho_{y_0,t_0}d\mu_t
	-\int_{\Sigma_t}Q\rho_{y_0,t_0}(\bar{H}\cdot(\bar{H}+E))d\mu_t \\
&+\int_{\Sigma_t}Q(-\Delta^{\Sigma_t}\rho_{y_0,t_0}-\rho_{y_0,t_0}
	(\frac{|\bar{F}^\perp|^2}{4(t_0-t)^2}+\frac{\bar{F}^\perp\cdot\bar{H}}{t_0-t}+\frac{\bar{F}^\perp\cdot E}{2(t_0-t)})d\mu_t
\end{align*} where we also use $\frac{d}{dt}d\mu_t = -|H|^2 d\mu_t = -\bar{H}\cdot(\bar{H}+E)d\mu_t$.

Integrating by parts, collecting terms and completing squares, we obtain
\begin{align*}
&\frac{d}{dt}\int_{\Sigma_t}Q\rho_{y_0,t_0}d\mu_t \\=&\int_{\Sigma_t}\big(P\rho_{y_0,t_0}-Q\rho_{y_0,t_0}(\frac{|\bar{F}^\perp|^2}{4(t_0-t)^2}
	+\frac{\bar{F}^\perp\cdot\bar{H}}{t_0-t}+\frac{\bar{F}^\perp\cdot E}{2(t_0-t)}+|\bar{H}|^2+\bar{H}\cdot E)\big)d\mu_t \\
=&\int_{\Sigma_t}\big(P\rho_{y_0,t_0}-Q\rho_{y_0,t_0}
	(|\bar{H}+\frac{\bar{F}^\perp}{2(t_0-t)}+E|^2-|E|^2)\big)d\mu_t.
\end{align*}

Recall (\ref{second_ineq}), $0\leq Q\leq 1$, and  that $\lim_{t\rightarrow t_0} \int_{\Sigma_t}\rho_{y_0,t_0}d\mu_t$ exists, we arrive at
$$ \frac{d}{dt}\int_{\Sigma_t}(1-\eta)\rho_{y_0,t_0}d\mu_t \leq
C-\frac{2}{3}\int_{\Sigma_t}\eta|\mbox{II}|^2\rho_{y_0,t_0}d\mu_t $$ for some constant $C>0$.
Moreover, for $t<\infty$, we know that $\eta>\delta>0$, thus
\begin{align}\frac{d}{dt}\int_{\Sigma_t}(1-\eta)\rho_{y_0,t_0}d\mu_t \leq
C-C_\delta\int_{\Sigma_t}|\mbox{II}|^2\rho_{y_0,t_0}d\mu_t \;\;\;\;\;\; \forall t \in [0,t_0) \label{monotoneineq}
\end{align}
From here we deduce that
$\lim\limits_{t\rightarrow t_0} \int_{\Sigma_t}(1-\eta)\rho_{y_0,t_0}d\mu_t$ exists.

Consider the parabolic dilation at $(y_0,t_0)$:
$$ \begin{array}{cccc}
D_\lambda: &\mathbb{R}^N\times[0,t_0) &\rightarrow &\mathbb{R}^N\times[-\lambda^2t_0,0) \\
&(y,t) &\mapsto &(\lambda(y-y_0),\lambda^2(t-t_0))
\end{array}. $$
Notice that $\eta$ is invariant under $D_\lambda$. The inequality (\ref{monotoneineq}) on $\Sigma^\lambda_s$ becomes
\begin{align*}
\frac{d}{ds}\int_{\Sigma^\lambda_s}(1-\eta)\rho_{0,0}d\mu^\lambda_s
\leq \frac{C}{\lambda^2}-C_{\delta}\int_{\Sigma^\lambda_s}|\mbox{II}|^2\rho_{0,0}d\mu^\lambda_s.
\end{align*}
Fix $s_0<0, \tau>0$, for $\lambda$ large, we integrate both sides over the interval $[s_0-\tau, s_0]$ and obtain
\begin{align*}
&\int_{s_0-\tau}^{s_0}\int_{\Sigma^\lambda_s}|\mbox{II}|^2\rho_{0,0}d\mu^\lambda_s ds
\leq \frac{C'}{\lambda^2}
+C''\int_{\Sigma^\lambda_{s_0-\tau}}(1-\eta)\rho_{0,0}d\mu^\lambda_{s_0-\tau}
-C''\int_{\Sigma^\lambda_{s_0}}(1-\eta)\rho_{0,0}d\mu^\lambda_{s_0}.
\end{align*}
 Since the limit $\lim_{t\rightarrow t_0}\int_{\Sigma_t}(1-\eta)\rho_{y_0,t_0}d\mu_t$ exists, argue as before we can pick a sequence $\lambda_i\rightarrow\infty$ such that
$\Sigma^{\lambda_i}_s\rightarrow\Sigma^\infty_s$ for all $s\in(-\infty,0)$. In fact, we have
$$ \int_{s_0-\tau}^{s_0}\int_{\Sigma^{\lambda_i}_s}|\mbox{II}|^2\rho_{0,0}d\mu^{\lambda_i}_s ds \leq C(i) $$
where $C(i)\rightarrow0$ as $i\rightarrow\infty$. We first choose $\tau_i\rightarrow0$ such that
$\frac{C(i)}{\tau_i}\rightarrow0$, and then choose $s_i\in[s_0-\tau_i,s_0]$ so that
$$ \int_{\Sigma^{\lambda_i}_{s_i}}|\mbox{II}|^2\rho_{0,0}d\mu^{\lambda_i}_{s_i}\leq\frac{C(i)}{\tau_i}
	\rightarrow 0 $$

Suppose $\Sigma^{\lambda_i}_{s_i}$ is given by $\bar{F}^{\lambda_i}_{s_i}:\Sigma \rightarrow \R^N$, thus
$$ \rho_{0,0}(F^{\lambda_i}_{s_i})=\frac{1}{4\pi(-s_i)}\exp\big(-\frac{|\bar{F}^{\lambda_i}_{s_i}|^2}{4(-s_i)}\big). $$
Then for any $R>0$,
$$ \int_{\Sigma^{\lambda_i}_{s_i}}\rho_{0,0}|\mbox{II}|^2d\mu^{\lambda^i}_{s_i}
\geq \int_{\Sigma^{\lambda_i}_{s_i}\bigcap B_R(0)}\rho_{0,0}|\mbox{II}|^2d\mu^{\lambda^i}_{s_i}
\geq C\exp(-\frac{R^2}{2})\int_{\Sigma^{\lambda_i}_{s_i}\bigcap B_R(0)}|\mbox{II}|^2d\mu^{\lambda^i}_{s_i} $$
Hence, on any compact set $K\subset\mathbb{R}^N$,
$$ \lim_{i\rightarrow\infty}\int_{\Sigma^{\lambda_i}_{s_i}\bigcap K}|\mbox{II}|^2d\mu^{\lambda^i}_{s_i}=0 $$

We can take a coordinate neighborhood $\Omega$ of $\pi_1(y_0)\in \Sigma^1$ where $\pi_1:M\rightarrow \Sigma^1$ is the projection onto the first factor. $\Sigma^{\lambda_i}_{s_i}$ is the
graph of $\tilde{u}_i:\lambda_i \Omega \rightarrow \lambda_i \Sigma_2$. Since $\eta$ is bounded and
$\int_{\Sigma^{\lambda_i}_{s_i}} |\mbox{II}|^2 d\mu^{\lambda_i}_{s_i}\rightarrow 0$,
$$ |D\tilde{u}_i|\leq C \mbox{ and } \int_\Omega |D^2\tilde{u}_i|^2\rightarrow0 $$
Therefore $\tilde{u}_i\rightarrow\tilde{u}_\infty$ in $C^\alpha\bigcap W^{1,2}$, $\tilde{u}_\infty$ is a linear
map, and
$$ \lim_{i\rightarrow\infty}\int_{\Sigma^{\lambda_i}_{s_i}} \rho_{0,0}d\mu^{\lambda_i}_{s_i}=\int\rho_{0,0}d\mu^\infty_{s_0}=1 $$
We thus found a sequence such that $\lim_{t_i\rightarrow t_0}\int\rho_{y_0,t_0}d\mu_t$ exists and
equals to $1$. By White's regularity theorem \cite{wh3}, the second fundamental form $|\mbox{II}|$ is bounded as $t\rightarrow t_0$.

\subsection{Smooth convergence as $t\rightarrow\infty$}
 By the general convergence theorem of Simon \cite{si}, it suffices to have an uniform bound on $\sup_{\Sigma_t}|\mbox{II}|^2$ as $t\rightarrow \infty$.

We first prove the following differential inequality.
\begin{lem}
$$ \frac{d}{dt}\int_{\Sigma_t}\frac{|H|^2}{\eta} d\mu_t\leq c\int_{\Sigma_t}\frac{|H|^2}{\eta}d\mu_t $$
\end{lem}
\begin{proof}
It is not hard to derive:
\[(\frac{d}{dt}-\Delta)|H|^2=-2|\nabla H|^2+2\sum_{i,j}(\sum_\alpha H_\alpha h_{\alpha ij} )^2+c(2-\eta^2)|H|^2\] where $H=H_\alpha e_\alpha$ for an orthonormal basis $e_\alpha$ of the normal space of $\Sigma_t$.
Combining this equation with (\ref{eta_lg}) and applying the inequalities $\sum_{i,j}(\sum_\alpha H_\alpha h_{\alpha ij} )^2\leq |\mbox{II}|^2|H|^2$ and $|\nabla|H||\leq|\nabla H|$, the inequality follows.
\end{proof}

The following proposition, which follows from a standard point-picking lemma, will be useful in the proof too.
\begin{prop}\label{blow_up_infty}
Suppose $\sup_{\Sigma_t}|\mbox{II}|^2$ is not bounded as $t\rightarrow \infty$, then there exists a blow-up flow
$\widetilde{S}_\infty\subset\mathbb{R}^4\times\mathbb{R}$ defined on the whole $(-\infty,\infty)$ with uniformly
bounded second fundamental form and $|\widetilde{\mbox{II}}|(0,0)=1$. Indeed, each $t$ slice of $\widetilde{S}_\infty$ is a graph of a symplectomorphism from $\mathbb{C}$ to $\mathbb{C}$.
\end{prop}

\subsection{$c>0$ case}
We prove by contradiction and look at the blow-up flow $\widetilde{S}_\infty$.
(\ref{eta_lower}) implies $\eta\equiv1$ on
$\widetilde{S}_\infty$. The $t=0$ slice of $\widetilde{S}_\infty$ is a graph of a symplectomorphism from $\C$ to $\C$ given by
\[(x,y)\mapsto (f(x,y),g(x,y))\] and satisfies
\[f_x g_y-f_y g_x=1 \text{  and   }
\frac{2}{\sqrt{2+f^2_x+f^2_y+g^2_x+g^2_y}}=1.\]
We derive that $h=f+\sqrt{-1}g$ must be holomorphic with $|\frac{\partial h}{\partial z}|=1$, $z=x+\sqrt{-1}y$. Thus $h$ is of the form
$h=e^{\sqrt{-1}\theta}z+C$  where $\theta$ and $C$ are constants. The graph of $h$ has zero second fundamental form
and this contradicts with $|\tilde{\mbox{II}}|(0,0)=1$.

\subsection{$c=0$ case}
 In this case by (\ref{eta_lower}), $\eta$ has a positive lower bound. We have
\begin{equation}\label{H}\int_{\Sigma_t} |H|^2 d\mu_t\leq
\int_{\Sigma_t}\frac{|H|^2}{\eta}d\mu_t \leq
C\int_{\Sigma_t}|H|^2 d\mu_t\end{equation} for some constant $C$.

Since $\int_0^\infty \int_{\Sigma_t} |H|^2 d\mu_t dt<\infty$, there exists
a subsequence $t_i$ such that $\int_{\Sigma_{t_i}} |H|^2 d\mu_{t_i}
\rightarrow 0$ and thus $\int_{\Sigma_{t_i}} \frac{|H|^2}{\eta}d\mu_{t_i}
\rightarrow 0$ as well. Because $\int_{\Sigma_t}
\frac{|H|^2}{\eta} d\mu_t$ is non-increasing, this implies
$\int_{\Sigma_{t}} \frac{|H|^2}{\eta}d\mu_t \rightarrow 0$ for the
continuous parameter $t$ as it approaches $\infty$.  Together with
(\ref{H}) this implies $\int_{\Sigma_{t}}{|H|^2} d\mu_t \rightarrow 0$ as
$t\rightarrow \infty$. By the Gauss-Bonnet theorem,
$\int_{\Sigma_{t}}{|\mbox{II}|^2}d\mu_t=\int_{\Sigma_t}|H|^2 d\mu_t\rightarrow 0$. The
$\epsilon$ regularity theorem in \cite{il1} (see also \cite{ec})
implies $\sup_{\Sigma_t} |\mbox{II}|^2 $ is uniformly bounded.

\subsection{$c<0$ case}
We may assume $c=-1$. In this case we have
\[\frac{d}{dt}\int_{\Sigma_t} \frac{|H|^2}{\eta}d\mu_t\leq -\int_{\Sigma_t}\frac{|H|^2}{\eta} d\mu_t
\]
or
\[\int_{\Sigma_t} \frac{|H|^2}{\eta}d\mu_t\leq C e^{-t}\] for some
constant $C$.

 Since $\eta\leq 1$, we have
\[\int_{\Sigma_t} {|H|^2}d\mu_t\leq C e^{-t}.\]

Suppose the second fundamental form is unbounded.  Since $\int_{\Sigma_t} |H|^2d\mu_t \leq C e^{-t}$,
the integral of $|H|^2$ on each time-slice of the blow-up flow $\widetilde{S}_\infty$ vanishes.
Therefore, we obtain a minimal area-preserving map. A result of Ni
\cite{ni} shows this is a
linear diffeomorphism. This contradicts to the fact that
$|\tilde{\mbox{II}}|(0,0)=1$ again.
\section{Acknowledgement}
These notes were prepared for a series of lectures the author gave at the Taipei branch of National Center of Theoretical Science in the summer of 2005. The author would like to thank Professor Ying-Ing Lee for her kind invitation, the center for the hospitality during his visit, and Mr. Kuo-Wei Lee and Mr. Chung-Jun Tsai for taking excellent notes and  typing a preliminary version of these notes.
The material is based upon work supported by the National Science
Foundation under Grant No. 0306049 and 0605115.


\begin{thebibliography}{aa}
\bibitem{ec} K. Ecker, \textit{Regularity theory for mean curvature
flow}, Progress in Nonlinear Differential Equations and their
Applications, 57. Birkhäuser Boston, Inc., Boston, MA, 2004.

\bibitem{eh1} K. Ecker and G. Huisken,
\textit{Mean curvature evolution of entire
graphs}, Ann. of Math. (2) \textbf{130} (1989)
, no. 3, 453--471.


\bibitem{eh2} K. Ecker and G. Huisken,
\textit{Interior estimates for hypersurfaces moving by mean
curvature}, Invent. Math. \textbf{105} (1991), no. 3, 547--569.


\bibitem{hu1} G. Huisken, \textit{Flow by mean
 curvature of convex surfaces into spheres},
J. Differential Geom. \textbf{20} (1984), no. 1, 237--266.


\bibitem{hu2} G. Huisken, \textit{Asymptotic
behavior for singularities of the mean curvature
flow}, J. Differential Geom. \textbf{31} (1990), no.
1, 285--299.

\bibitem{il1} T. Ilmanen, \textit{Singularities
of mean curvature flow of surfaces},
preprint , 1997.

\bibitem{il2} T. Ilmanen, \textit
{Elliptic regularization and partial
regularity for motion by mean curvature},
 Mem. Amer. Math. Soc. \textbf{108} (1994), no. 520,


\bibitem{le} Y.-I. Lee, \textit{Lagrangian minimal surfaces in K\"ahler-Einstein surfaces of negative scalar curvature}, Comm. Anal. Geom. \textbf{2} (1994), no. 4, 579--592.


\bibitem{ni} L.  Ni, \textit{A Bernstein type theorem for minimal volume preserving
maps}, Proc. Amer. Math. Soc. \textbf{130} (2002), no. 4, 1207--1210.

\bibitem{sc} R. Schoen, \textit{The role of harmonic mappings in rigidity and
 deformation problems}, Complex geometry
 (Osaka, 1990), 179--200, Lecture Notes in
 Pure and Appl. Math., 143, Dekker, New York,
 1993.







\bibitem{si} L. Simon, \textit{Asymptotics
 for a class of nonlinear evolution
 equations, with applications to geometric
  problems}, Ann. of Math. (2) \textbf{118} (1983),
  no. 3, 525--571
\bibitem{sm1} K. Smoczyk,
\textit{A canonical way to deform a Lagrangian submanifold},
preprint, dg-ga/9605005.

\bibitem{sm2} K. Smoczyk, \textit{Angle theorems for the Lagrangian mean curvature
flow.} Math. Z. \textbf{240} (2002), no. 4, 849--883.

\bibitem{wa1} M.-T. Wang, \textit{Mean curvature flow of surfaces
 in Einstein four-manifolds}, J.Diff.Geom. {\textbf 57}(2001), no. 2, 301-338.
\bibitem{wa2} M.-T. Wang, \textit{Deforming area preserving
diffeomorphism of surfaces by mean curvature flow}, Math. Res. Lett. \textbf{8} (2001),
no.5-6, 651-662.


\bibitem{wa3} M.-T. Wang, \textit{A convergence result of the Lagrangian mean curvature
flow}, in the Proceedings of the third
 International Congress of Chinese Mathematicians. arXiv:math/0508354.
\bibitem{wh1} B. White, \textit{The size of the
singular set in mean curvature flow of mean convex
surfaces}, J. Amer. Math. Soc. \textbf{13} (2000), no. 3, 665--695.
\bibitem{wh2} B. White, \textit{Stratification
of minimal surfaces, mean curvature flows, and
 harmonic maps}, J. Reine Angew. Math. \textbf{488}
 (1997), 1--35.

\bibitem{wh3} B. White, \textit{A local
regularity theorem for classical mean curvature
flow}, Ann. of Math. (2) \textbf{161} (2005), no. 3, 1487--1519.


\end{thebibliography}
\end{document}